\documentclass[12pt]{amsart}
\usepackage{amsmath}
\usepackage{amsthm}
\usepackage{amsfonts}

\theoremstyle{remark}

\def\intzc #1'{{\text{$int_{#1} A^c$}}}
\def\rat{\mathbb{Q}}
\def\real{\mathbb{R}}
\def\rmoc #1{{\text{$\real^#1$}}}
\def\m #1{{\text{$H_#1$}}}

\def\dhz {{\text{$\overline{d} (x,A) \ $}}}

\def\ddz {{\text{$\underline{d} (x,A) \ $}}}
\def\fr #1'{{\text{$fr_{#1}A$}}}
\def\intp #1,#2'{{\text{$int_{#1}{#2}$}}}
\def\intz #1'{{\text{$int_{#1} A$}}}
\def\je {{\text{$\Im_1^{n-1}$}}}
\def\fa #1{{\text{$\Phi_{\Omega,#1}^g$}}}

\def\pez #1'{{\text{$ P_{\Omega,#1} $}}}
\def\pe {{\text{$ \pez \tau' $}}}

\def\pda {{\text{$\pe (A) $}}}
\def\drr #1'{{\text{$\Omega \cap \fr #1' $}}}
\def\mitau #1'{{\text{$\mu_\tau(\drr #1')$}}}
\def\mien #1'{{\text{$\mu_n(\drr #1')$}}}
\def\en  {{\text{$N(p_\tau,\drr e',z)$}}}
\def\enpr {{\text{$N(p_n,\drr pr',z)$}}}
\def\mad {{\text{$m^A_{\Omega,\tau}$}}}
\def\zsz {{\text{$z+\frac sm e_n$}}}
\def\zsm #1,#2'{{\text{$z_{#1} + \frac{s+#2}{m}e_{n}$}}}
\def\ins {{\text{$\left]\frac{s-1}{m},\frac{s+1}{m}\right]$}}}
\def\ti #1'{{\text{$ \tau_{#1} $}}}
\def\rn {{\text{$r_0$}}}
\def\uxr #1'{{\text{$U(x_{#1},\rn)$}}}
\def\uxrn {{\text{$\uxr 0'$}}}
\def\rpodm #1{{\text{$#1 \subset \rmoc n$}}}
\def\rnt {{\text{$\real^{n-1}(\tau)$}}}
\def\ip {{\text{$int_{pr}A$}}}
\def\madn {{\text{$m^{\ip}_{\Omega,n}$}}}
\def\apl {{\text{$A^+(k,m,s)$}}}
\begin{document}
	
	\author{Miroslav Chleb\'\i k\\ \ \\}
	
	\address{University of Sussex, UK} 
	\email{m.chlebik@sussex.ac.uk}

	\title{GOING BEYOND VARIATION OF SETS}

\begin{abstract} We study integralgeometric representations of variations 
of general sets $A \subset \rmoc n$ without
any regularity assumptions. If we assume, for example,
that just one partial derivative of its characteristic function 
$\chi^A$ is a signed Borel measure on $\rmoc n $
with finite total variation, can we provide
a nice integralgeometric representation of this variation?
This is a delicate question, as the Gauss-Green type theorems
of De Giorgi and Federer are not  available in this generality.
We will show that a `measure-theoretic boundary' plays its role 
in such representations similarly as for the sets of finite variation.
There is a variety of suitable notions of  `measure-theoretic boundary'
and one can address the question to find notions of measure-theoretic boundary 
that are as fine as possible.

\end{abstract}

\keywords{perimeter of sets, measure-theoretic boundary, integralgeometric measure}

\subjclass{Primary: 28A75, 49Q15; Secondary: 26B15, 28A78}

\maketitle

\section{INTRODUCTION}

 A recurring theme in Geometric Measure Theory and in the study of geometric variational problems is the theory of sets of finite perimeter. The best known classical result about such sets due to De Giorgi and Federer (\cite{DG1}, \cite{DG2}, \cite{FH1}, and \cite{FH2})
says
that a (measurable) set $A \subset \rmoc n$ has finite perimeter if and only if its measure-theoretic boundary has finite area ($(n-1)$-dimensional Hausdorff measure), and more precisely the perimeter agrees with the area of the measure-theoretic boundary of $A$.

  Our focus here is on more general framework of sets. The main result in this paper (Theorem 4.5) states that a set $A$ has finite variation in a given direction $\tau$ (that is, the distributional derivative of the characteristic function of $A$ in the direction $\tau$ is a finite measure) if and only if a suitably defined $(n-1)$-dimensional measure of a suitably defined measure-theoretic boundary is finite, and more precisely the variation of $A$ in  the direction $\tau$ agrees with the measure of such boundary. Interestingly, our results give also a relatively elementary proof of the classical result of De Giorgi and Federer mentioned above (Theorem 4.9) . The results show quite clearly that the natural notion of area is not the $(n-1)$-dimensional Hausdorff measure, but the integralgeometric measure (which of course agree in case of rectifiable sets).

A set $ A \subset \rmoc n $ is said to be of finite
perimeter if it is Lebesgue measurable and 
the gradient $D\chi^A$ in the sense of distributions 
of its characteristic function $\chi^A$ is an $\rmoc n$
valued Borel measure on $\rmoc n $ with finite total variation.
The value of the perimeter of $A$, denoted by $P(A)$, is then
the total variation $||D\chi^A||$ of the vector measure $D\chi^A$.
Otherwise, let the perimeter of $A$ be equal to  + $\infty$. 
(Another equivalent
definition of perimeter was given in  \cite{DG1}, see also \cite{DG2}
and \cite{FH2}.)

Given a direction $\tau \in S^{n-1}$ a set $ A \subset \rmoc n $ is said 
to have bounded variation at the direction $\tau$ if it is Lebesgue measurable
and the directional derivative  in the sense of distributions $\partial_\tau \chi^A$
of its characteristic function $\chi^A$ is a signed Borel measure with finite total variation
on 
$\rmoc n $ . The value of the variation at direction $\tau$ of $A$,
denoted by $P_\tau (A)$, is then the total variation $||\partial_\tau \chi^A||$ of 
the signed measure $\partial_\tau \chi^A$. Otherwise, let  $P_\tau (A) = + \infty$.

It is well known that, for a Lebesgue measurable set  $A$ and $\tau = e_{i}$
being the standard orthonormal basis direction (and writing briefly $P_{i}$
instead of  $P_{e_i}$),

$$ P_{i} (A) = \int m_i^A(z) \, dz $$
where $ m_i^A (z)$  is the infimum of the variations in $x_i$ of all
functions defined on the line $ L_i (z)$  (parallel to the $x_i$
axis and meeting $z$)  which are equivalent to $\chi^A \vert L_i(z) $
and the integration is over the $(n-1)$ space orthogonal to
the $x_i$ axis.

  It is known that the perimeter of $A$ (if it is finite) is
equal to the $(n-1)$ measure  of the set \fr r'  that is called the
reduced boundary (see \cite{FH2}) or equivalently it is equal
to $(n-1)$ measure of the essential boundary \fr e' of $A$ (see
\cite{VA} or \cite[4.5.6]{FH3}). Specifically, $x \in \fr r'$ iff there is
an $(n-1)$ plane $\pi$ through $x$  such that the symmetric
difference of $A$ and one of the halfspaces determined by $\pi$
has density zero at $x$. Further, $x \in \fr e'$  iff both $A$ and
complement of $A$ have positive outer upper density at $x$.

  Moreover, if the $(n-1)$ measure 
 of \fr e' is finite then $A$ is of finite perimeter \cite[4.5.11]{FH3}.
Hence $(n-1)$
measure of  \fr e'  is equal to the perimeter of $A$ for a general
set $ A \subset \rmoc n$ (Our method also
offers a simple self-contained proof of this fact for
an integralgeometric $(n-1)$ measure.)

The main purpose of this paper is to show that the directional
variation of a general set $ A \subset \rmoc n$ (without any
assumptions on regularity of $A$) is equal to the measure
of projection (with multiplicities taken into account) of 
the 'measure-theoretic boundary'. The essential boundary \fr e'
can play here the role of such a 'measure-theoretic
boundary', but one can aim to replace it even with finer notions of
'measure-theoretic boundary'.
We show that one can replace  \fr e'  by finer preponderant
boundary  \fr pr'  (see 4.5, 4.8, 4.9 and 4.10). Specifically,
$x \in \fr pr'$ iff both $A$ and complement of $A$ have the outer
upper density at $x$ greater than or equal to $\frac 12$.

\section{NOTATION AND TERMINOLOGY}

Throughout the whole paper we
deal with the sets in the n-dimensional Euclidean space $ \real^{n} $.
We tacitly assume that  $n \geq 2 $  but results
trivially  hold  in the case  $n = 1 $.

  Let  $ e_{1} , e_{2} , \ldots , e_{n} $ stand for the orthonormal
base in \rmoc n,
$
   e_{1}  = \left( 1 , 0 , 0 , \ldots , 0 \right)$,
  $ e_{2}  = \left( 0 , 1 , 0 , \ldots , 0 \right)
   , \ldots $, 
$   e_{n}  = \left( 0 , 0 , 0 , \ldots , 1 \right)$.
The symbol $\bold{0}$ has also the meaning of the zero vector $ \left(0,0,\ldots,0\right)
\text{ in } \real^{n} $.  For   $x,y \in \real^{n} $ we  denote  by  $ \vert x \vert $
the euclidean norm of  $x$  and by  $ x \circ y $  the inner product
of $x$ and $y$. The symbol $ \left[x,y\right]$  stands for the convex hull of
the set $ \left\{ x,y \right\} $ and $ \left] x,y \right[$  means $ [x,y] \setminus \{x,y\}$.

  Whenever  $x \in \rmoc n  $ and $ r > 0 \quad B(x,r)$  and  $U(x,r)$  stand
for the closed and open balls, respectively, with center $x$
and radius $r$ and  $Q(x,r) $ stands for the cubic interval
$$ \{ \ y \in \rmoc n :\ \left\vert y_i - x_i \right\vert \leq r \text{ , } 1 \leq i \leq n \  \} \, . $$
We put
$$
   S^{n-1}  =  \{ \ x \in \rmoc n  :\  \left\vert x \right\vert = 1  \ \} \quad \text{and} \quad
   L_\tau (x)  =   \{ \ x + t\tau  :\  t \in \real \ \} \ \ \text{for} \  x \in \rmoc n  \  \text{and} \  \tau \in %
   \rmoc n  \setminus \{ \bold{0} \} \,.
$$
For $\tau \in \rmoc n  \setminus \{ \bold{0} \} $ we denote by $\real^{n-1} (\tau)$
the orthogonal complement  in  $\rmoc n$ to the one  dimensional subspace
$ \{ t\tau : t \in \real \} $
and by  $p_\tau$  the orthogonal projection of $\rmoc n$ onto $\real^{n-1} (\tau)$.
We  write  briefly  $L_i(x)$, $ \real^{n-1} (i)$ and $p_i$ in the
case  $ \tau = e_i$.

 For $ A, B \subset \rmoc n $ we  denote  by $ A \bigtriangleup B $ the symmetric
difference of $A$ and $B$,
$   A \bigtriangleup B = ( A \setminus B ) \cup ( B \setminus A), $
and by $ A^c $ and  $\chi^A$ the complement of $A$ to  $\rmoc n $ and the
characteristic function of $A$ (on $\rmoc n$), respectively.

  For an open set  $\Omega \subset \rmoc n $ we will denote by  $ C^{\infty}_0 (\Omega)$  (and
 $ C^{\infty}_0 (\Omega,\rmoc n )$)
the space of all infinitely differentiable real
valued functions  with compact support in $\Omega$ (and the space
of all infinitely differentiable $\rmoc n$ valued vector
functions with compact support in $\Omega$, respectively).
These spaces are considered to be equipped with the ``sup norm''.

  For any function $f$, any set $A$ and any value $y$, the
multiplicity $ N(f,A,y) $ is defined as the number of
elements  (possibly  $+\infty$) of the set  $\{ x \in A : f(x) = y \}$.

  For $ A \subset \rmoc n  $ we denote by $int A , cl A \text{ and } fr A $
the interior, closure and boundary of $A$, respectively.\\

For any outer measure $\mu$ on $\rmoc n $ and for any
set  $X \subset \rmoc n $ we define the outer measure $\mu {\mid} X  \text{ on } \rmoc n $
by the formula
$$ ( \mu {\mid} X ) (A) = \mu \, ( X \cap A ) \quad \text{for every} \  A \subset \rmoc n. $$

For a signed Borel measure $\mu $ (or for a vector Borel
measure  $\mu = ( \mu_1,\mu_2,\ldots,\mu_n )$)  on \rmoc n the symbol $\left\vert \mu \right\vert$
stands for Borel measure which is the variation  of $\mu$ (or the variation of $\mu$ in the sense of vector measures, respectively).\\

{\bf 2.1 Hausdorff measures.} For an integer $k = 0,1,\ldots,n$ let
$\m k$  stand for the  $k$-dimensional Hausdorff outer measure on \rmoc n, 
which is normalized in such a way that
$$
 \m k \{ \ x \in \rmoc n:\  0 \leq x_i \leq 1 \quad \text{for} \quad 1 \leq i \leq k \quad \text{and }
                              x_i = 0 \quad \text{for} \quad k < i \leq n \ \} = 1.
$$
In particular, \m 0 is the counting measure and
\m n coincides with the Lebesgue outer measure on \rmoc n. 

  The constant $ V (n) =  \frac{\pi^\frac n2}{\Gamma(\frac n2 +1)} $ means the volume of the
unit ball in  \rmoc n (with $ V (0) = 1$),
and the constant
$ A(n) = nV(n) = \frac{2\pi^\frac n2}{\Gamma(\frac n2)} \
         \text{means the area of}  \ S^{n-1}  . $

 We define the equivalence  relation $\sim$ for subsets of  \rmoc n by the
prescription
$$
 A \sim B \quad \text{iff} \quad \m n \left[ A \bigtriangleup B \right] = 0  .
 $$

{\bf 2.2 Projection measures $\mu_\tau$.}\
 For  $\tau \in \rmoc n \setminus \{ \bold{0} \} $  the result
of  Caratheodory's  construction  from  the  set  function
$$ B \longmapsto  H_{n-1}[ \ p_\tau(B) \ ] $$
which is defined on the covering family of all Borel sets in
\rmoc n will be called  the projection measure at the direction
$\tau$ and denoted by $\mu_\tau$. Then  $\mu_\tau$ is a Borel regular
outer measure on \rmoc n and  $\mu_\tau \leq H_{n-1}$.

From Fubini theorem it follows that $H_n(C)=0$ whenever
$C \subset \rmoc n $ is such that $\mu_\tau(C)<\infty.$\\

{\bf 2.3 Integralgeometric measure \je.}\
 The  result  of
Caratheodory's   construction  from  the  set  function
 $$ B \longrightarrow \frac{1}{2V(n-1)} \int\limits_{S^{n-1}} H_{n-1}[ \ p_\tau(B) \ ] \, d H_{n-1}(\tau) $$
which is defined on the covering family of all Borel sets in
\rmoc n is usually termed  $(n-1)$ dimensional integralgeometric
measure with exponent 1 on  \rmoc n and denoted by \je . (For the
existence of the above integral see e.g.\  \cite[2.10.5]{FH3}.)

\je is a Borel regular outer measure on \rmoc n  and
$2V(n-1)\je \leq A(n)H_{n-1}$. Moreover $\je \leq H_{n-1}$ by \cite[3.3.16]{FH3}.\\

{\bf 2.4 Densities.}\
  For every set $A \subset \rmoc n $ and each  $x \in \rmoc n $
  we define the upper outer density \dhz and the lower
  outer density \ddz of $A$ at $x$ by the formulas
  \begin{align*}
    \dhz & = \overline{\lim}_{r \to 0+} \frac{\m n [A \ \cap \ B(x,r)]}{\m n [B(x,r)]} \ , \\
    \ddz & = \underline{\lim}_{r \to 0+} \frac{\m n [A \ \cap \ B(x,r)]}{\m n [B(x,r)]}.
 \end{align*}
 In the case  $\dhz = \ddz $ this common value is termed
 the outer density of $A$ at $x$ and it is denoted  by  $d(x,A)$.

  A point $x$ for which $\ddz = 1$ is termed the outer  density  point  of~$A$.
(We may drop the adjective ``outer'' from this terminology whenever the set
$A$ is Lebesgue measurable.)\\

{\bf 2.5 Essential and preponderant interior and boundary.}\
 We define the essential interior  $int_e A$ and the essential
boundary \fr e' of~the set  $A \subset \rmoc n $ by the formulas
$$\intz e' = \{ \ x \in \rmoc n :\ d(x, A^c) = 0 \},$$
$$\fr e'  =  \left[ \ \intz e'  \cup \intp e, A^c' \ \right]^c =
           \{ \ x \in \rmoc n :\ \dhz > 0 \quad \text{and} \quad \overline{d}(x, A^c) > 0 \ \}.$$
 It is easy to see that
 $\intz e'  \cap \intp e,(A^c)' = \emptyset,$
$\intz e' $ is of type $F_{\sigma \delta}$ and  \fr e'  is of type $G_{\sigma \delta} \ .$
We also define the preponderant interior \intz pr' and the
preponderant boundary \fr pr' of $A \subset \rmoc n $ by the formulas
$$\intz pr' =  \left\{ x \in \rmoc n : \overline{d}(x, A^c) < \frac 12 \right\},$$
$$\fr pr'=  \left[ \ \intz pr' \cup \intp pr,A^c' \ \right]^c=
\left\{ x \in \rmoc n : \dhz \geq \frac 12 \quad \text{and} \quad \overline{d}(x, A^c) \geq \frac{1}{2} \right\}.$$

 It is easy to see that $\intz pr' \cap \intp pr,A^c' = \emptyset,$
$\intz pr' $ is of type $F_{\sigma}$ and \fr pr' is of type $G_{\delta}$.\\
 
{\bf 2.6 BV  functions.}\
For a nonempty open set $ \Omega \subset \rmoc n $ and for any
$\tau \in \rmoc n $ we define the space $ BV(\Omega,\tau)$ of all locally (in $\Omega$)
\m n summable  functions  $g$  for  which  there exists a finite
signed Borel measure $ \Phi_{\Omega,\tau}^g$ on $\Omega$ with the equality
$$
 \int\limits_\Omega g(x) \cdot \tau \circ \operatorname{grad} \varphi (x) \, dx =
  - \int\limits_\Omega \varphi (x) \, d\Phi_{\Omega,\tau}^g (x)
$$
whenever $ \varphi \in C^{\infty}_0 (\Omega)$.
$BV(\Omega)$ is defined as the space of all locally (in $\Omega$) \m n
summable functions $g$ such that there exist the finite signed
Borel measures  $\fa 1, \fa 2 , \ldots , \fa n $ with the equality
$$
 \int\limits_\Omega g(x) \cdot \operatorname{div} \psi (x) \, dx =
 -\sum\limits_{i=1}^n \int\limits_\Omega \psi_{i} (x) \, d\fa i (x)
$$
whenever $\psi =(\psi_1, \psi_2,\dots,\psi_n)
\in C^{\infty}_0 (\Omega,\rmoc n )$.

We also define the space $ BV_{loc}(\Omega)$ (and the space
$BV_{loc}(\Omega,\tau)$ analogously) by the following:
$ g \in BV_{loc}(\Omega)$  iff $ g \vert \tilde \Omega \in BV(\tilde \Omega)$  whenever $ \tilde \Omega$  is nonempty
open set the closure of which is a compact subset of $\Omega$\\

{\bf 2.7 Directional variation and perimeter of sets.}\
 For a nonempty open set  $\Omega \subset \rmoc n$ and for any
$\tau \in \rmoc n $ the set functions
\pe \ and  $P_{\Omega}$ over the subsets of \rmoc n are defined  for
$A \subset \rmoc n$ by the following:
\begin{itemize}
\item
 If $ A \cap \Omega$  is not \m n measurable then we put
  $$ \pda = P_\Omega(A) = \infty.$$
\item
If  $A \cap \Omega $ is  \m n  measurable then we put
\end{itemize}
\begin{align*}
\pda & = \operatorname{sup} \left\{ \int\limits_{\Omega} \chi^A (x) \tau \circ D\varphi (x)  \, dx :\ \varphi \in C^{\infty}_0 (\Omega) \quad \text{and} \quad \vert \varphi \vert \leq 1 \right\}  , \\
P_{\Omega}(A) & = \operatorname{sup} \left\{ \int\limits_{\Omega} \chi^A (x) \operatorname{div}%
      \psi(x) \, dx :\ \psi \in C^{\infty}_0 (\Omega,\rmoc n ) \quad \text{and} \quad \vert \psi \vert \leq 1 \right\}  .
\end{align*}

The value  \pda \ is termed the variation at direction $\tau$
of the set $A$ in $\Omega$, and $P_{\Omega}(A)$ is the perimeter
of $A$ in $\Omega$.

  In  the case $ \chi^A \vert \Omega \in BV(\Omega,\tau)$  the symbol $ \Phi_{\Omega,\tau}^A$ stand for the
(uniquely determined) signed Borel measure on $\Omega$ such that
$$
 \int\limits_\Omega \chi^A (x) \tau \circ D\varphi (x) \, dx =
 -\int\limits_\Omega \varphi (x) \, d \Phi_{\Omega,\tau}^A (x)
$$
holds whenever $\varphi \in C^{\infty}_0 (\Omega)$. We write briefly $ \Phi_{\Omega,i}^A$ in the case $\tau = e_i $.

\section{AUXILIARY RESULTS}

{\bf 3.1 Lebesgue outer density theorem.}\
\rm For any set $ A \subset \rmoc n  $ \m n
almost every point of $A$ is an outer density point of~$A$.\\

{\bf 3.2 Remark.}\ Let $ \Omega \subset \rmoc n $ nonempty  open  and  let
$A \subset \rmoc n  $ arbitrary. Using  Lebesgue density theorem and Borel
regularity of Lebesgue outer measure
one can easily show that the following statements are true: \newline
If $ A \cap \Omega $ is  \m n  measurable  then
$$
 \Omega \cap \intz e'   \sim    \Omega  \cap  \intz pr'  \sim  \Omega \cap A \quad \text{and} \quad
 \m n (\Omega \cap \fr e' )  =   \m n  (\Omega  \cap  \fr pr' )  =  0.
$$
If $ A \cap \Omega $ is not  \m n measurable  then
$$ \m n \{ \ x \in \Omega \: d(x, A) = 1 \quad \text{and} \quad d(x, A^c) = 1 \ \} > 0 $$
and especially
$$ \m n  ( \Omega \cap \fr e' )  \geq \m n ( \Omega \cap \fr pr' ) > 0 \, .$$\

{\bf 3.3 Observations.}\	
\begin{itemize}
\item[(1)]
  If $ g \in BV(\Omega) $ then $ g \in BV(\Omega,\tau ) $ for every $\tau \in \rmoc n \setminus \{ \bold{0} \} $.


\item[(2)]
  If $ \tau_j \in \rmoc n , \alpha_j \in \real , g \in BV(\Omega,\tau_j ) \ \  (j = 1,2, \dots , r ) $
   and $ \tau =  \sum^r_{i=1} \alpha_j \tau_j $ then $ g \in BV(\Omega,\tau ) $ and
   $\Phi^g_{\Omega,\tau} = \sum^r_{i=1} \alpha_j \Phi^g_{\Omega,\tau_j }$.
\item[(3)]
 If  $\tau_1 , \tau_2 , \dots , \tau_n  \in \rmoc n $ are linearly independent and
 $g \in BV(\Omega,\tau_j )$ ($j = 1, 2, \dots ,n $)  then $g \in BV(\Omega).$
\item[(4)]
 $\pda < \infty $ holds if and only if $\chi_A \vert \Omega \in BV(\Omega,\tau )$. In this case
 $\pda = \vert \Phi^A_{\Omega,\tau} \vert (\Omega) $ holds.
\item[(5)]
 $P_\Omega(A) < \infty $ holds if and only if $\chi_A \vert \Omega \in BV(\Omega)$. In this case
 $P_\Omega(A) = \vert \Phi^A_\Omega \vert (\Omega) $, where $\Phi^A_\Omega = (\Phi^A_{\Omega,1},\Phi^A_{\Omega,2}, \dots ,\Phi^A_{\Omega,n} ) $.
\item[(6)]
If $P_\Omega(A) = \infty $ then the set $\{ \tau \in \rmoc n:\ \pda < \infty \}$ is contained in an
$(n-1)$-dimensional linear subspace of \rmoc n.
 \end{itemize}

\bigskip

{\bf 3.4 Lemma.}\ Let  $B \subset \rmoc n $ be a Borel set.
	
\begin{itemize}
\item[(i)]
For any $\tau \in S^{n-1}$ the function $ z \longmapsto N(p_\tau,B,z) $
defined on  $\real^{n-1}(\tau)$ is $ H_{n-1}$ measurable and
 $$ \mu_\tau(B) = \int\limits_{\real^{n-1}(\tau)} N(p_\tau,B,z) \, dH_{n-1}(z).$$

\item[(ii)]
The function  $ \tau \longmapsto \mu_\tau(B) $ defined on $S^{n-1}$ is 
$ H_{n-1}$ measurable and
$$ \je (B) = \frac{1}{2V(n-1)} \int\limits_{S^{n-1}}\mu_\tau(B) \, dH_{n-1}(\tau)  . $$
\end{itemize}

\begin{proof}
See   \cite[2.10.10]{FH3} and  \cite[2.10.15]{FH3}.
\end{proof}

{\bf 3.5 Definition.}
\rm Let  $L_0$ be a line in \rmoc n and let $ L \subset L_0 $
be relatively open in $L_0$. For any set $ A \subset \rmoc n$ the
point $ x \in L $ is termed a hit of $L$ on $A$ provided both
$L \cap A \cap U(x,r)$ and  $( L \setminus A ) \cap U(x,r) $ have a positive
\m 1 measure for every $ r > 0 $. 

For $ \Omega \subset \rmoc n $ nonempty  open, $ A \subset \rmoc n $ , $z \in \rmoc n $
and $\tau \in S^{n-1} $ the symbol $ M_{\Omega,\tau}^A(z)$ stands for the set of all
hits of $ L_\tau (z) \cap \Omega$ on A $, m_{\Omega,\tau}^A(z)$  stands for the number of
elements (possibly $ +\infty$)  of $ M_{\Omega,\tau}^A(z)$ and  we put
$$  M_{\Omega,\tau}^A = \cup \, \{ \ M_{\Omega,\tau}^A(z) \: z \in \rmoc n \ \}.$$
We write briefly $ M_{\Omega,i}^A(z) $,  $m_{\Omega,i}^A(z)$  and $ M_{\Omega,i}^A$ in the case
$\tau = e_i$.\\

The starting point to our results is the following  known integral representation of directional variation of a set.\
\bigskip

{\bf 3.6 Lemma.}
\rm Let $ \Omega \subset \rmoc n $ be  nonempty  open  and  let
$ A \subset  \rmoc n$ be such that $ A \cap \Omega$ is \m n  measurable. Then for every
$\tau \in S^{n-1}$ the function $  z \longmapsto  m_{\Omega,\tau}^A(z)$ defined on
$\real^{n-1}(\tau)$ is $H_{n-1}$ measurable  and
$$ \pda =   \int\limits_{\real^{n-1}(\tau)}  m_{\Omega,\tau}^A(z)  \, dH_{n-1}(z).$$

\begin{proof}
See  e.g. \cite{MJ} and Chap.~7 of \cite{KK}.
\end{proof}

\section{INTEGRALGEOMETRIC CHARACTERIZATION OF VARIATIONS}

{\bf 4.1 Notation.}\ Let $ \Omega \subset \rmoc n  $ be  nonempty  open  and  $ A \subset \rmoc n  $
be such  that  $ A \cap \Omega $ is  \m n  measurable.  Let us identify  \rmoc n  with $\real^{n-1} \times \real$. For any $\alpha, \beta$ such that $-\infty \leq \alpha < \beta  \leq+ \infty$ put

$$
 E_{\Omega}(\alpha, \beta; A) = \{ \ z \in \real^{n-1}  \, : \,   \{z\} \times (\alpha, \beta) \subset \Omega \quad \text{and} \quad
                \m 1 ( \{z\} \times (\alpha, \beta) \setminus A ) = 0 \ \}  .
$$

It easily follows from Fubini's theorem that all the sets $E_{\Omega}(\alpha, \beta; A)$ in $\real^{n-1}$ are $H_{n-1}$ measurable.\\

{\bf 4.2 Lemma.}\
Let $ \Omega \subset \rmoc n  $ be  nonempty  open  and  $ A \subset \rmoc n  $
be such  that  $ A \cap \Omega $ is  \m n  measurable. Then there is an $H_{n-1}$ null set
$N \subset \real^{n-1}$ such that every $z \in \real^{n-1} \setminus N$ has the following properties:

\begin{itemize}
\item[(a)]\ If $\alpha, \beta \in \rat \cup \{-\infty, +\infty\}, -\infty \leq \alpha< \beta\leq+\infty$
($\rat$ being the set of rationals) are such that $z \in E_{\Omega}(\alpha, \beta; A)$ 
($z \in E_{\Omega}(\alpha, \beta; A^c)$, respectively) then $z$ is a density point in $\real^{n-1}$
of $E_{\Omega}(\alpha, \beta; A)$ (of $E_{\Omega}(\alpha, \beta; A^c)$, respectively).

\item[(b)]\  If $-\infty \leq \alpha< \beta\leq+\infty$ are such that $\{z\} \times (\alpha, \beta) \subset \Omega$
and $\m 1 ( \{z\} \times (\alpha, \beta) \setminus A ) = 0$ ($\m 1 ( \{z\} \times (\alpha, \beta) \setminus A^c ) = 0$, respectively) then $\{z\} \times (\alpha, \beta) \subset \intz e'$ ($\{z\} \times (\alpha, \beta) \subset \intzc e'$, respectively).\

\item[(c)]\   $\{x \in \Omega \cap  \fr e' :\ p_n(x)=z \} \subset M_{\Omega,n}^A.$
\end{itemize}

\begin{proof}

\begin{itemize}
\item[(a)]\  For any $H_{n-1}$ measurable set $B \subset \real^{n-1}$ put
$\tilde B=\{z \in B: z \ \text{is not a density point of} \ B \}.$ Due to the Lebesgue density theorem
$\tilde B$ is an $H_{n-1}$ null set. Hence the set
$$N= \cup \{ \tilde E_{\Omega}(\alpha, \beta; A) \cup  \tilde E_{\Omega}(\alpha, \beta; A^c):
\alpha, \beta \in \rat \cup \{-\infty,+\infty \}, \alpha<\beta \}$$

is an $H_{n-1}$ null set and each $z \in  \real^{n-1} \setminus N$ has the property (a).

\item[(b)]\ If $z \in  \real^{n-1} \setminus N$ and $-\infty \leq \alpha < \beta \leq +\infty$
are such that $\{z\} \times (\alpha, \beta) \subset \Omega$
and $\m 1 ( \{z\} \times (\alpha, \beta) \setminus A ) = 0$ ($\m 1 ( \{z\} \times (\alpha, \beta) \setminus A^c ) = 0$, respectively), then $z$ is a density point in $\real^{n-1}$
of $E_{\Omega}(\alpha_1, \beta_1; A)$ (of $E_{\Omega}(\alpha_1, \beta_1; A^c)$, respectively)
whenever $\alpha_1, \beta_1 \in \rat \cup \{-\infty, +\infty\}$ with
$\alpha \leq \alpha_1 < \beta_1\leq \beta$. From Fubini's theorem it follows that
$\{z\} \times (\alpha, \beta) \subset \intz e'$ ($\{z\} \times (\alpha, \beta) \subset \intzc e'$, respectively),
hence (b) holds true.

\item[(c)]\ Let us keep $z \in  \real^{n-1} \setminus N$ fixed and assume that 
$x \in \Omega \setminus M_{\Omega,n}^A$ is such that $p_n(x)=z$. Our aim is to prove
that then necessarily $x \notin \fr e'$. As $x \in \Omega \setminus M_{\Omega,n}^A$ we can
find real numbers $\alpha<\beta$ such that $x \in \{z\} \times (\alpha, \beta) \subset \Omega$
and either $\m 1 ( \{z\} \times (\alpha, \beta) \setminus A ) = 0$ or $\m 1 ( \{z\} \times (\alpha, \beta) \setminus A^c ) = 0$.  From (b) it follows that either $x \in \intz e'$ or $x \in \intzc e'$, hence 
$x \notin \fr e'$. This completes the proof.

\end{itemize}

\end{proof}

{\bf 4.3 Lemma.}\ Let $  X,Y \subset \real $ be two disjoint sets
of type $ F_{\sigma} $ such that every $ x \in X $ is a bilateral
accumulation point of $\real \setminus Y $ and every $ y \in Y $  is a bilateral
accumulation  point  of $\real \setminus X $. Then $  ] \, a,c \, [ \, \setminus \, ( X \cup Y ) $ is
nonempty whenever $ a \in X $ and  $ c \in Y $.

\begin{proof}
We have $ X = \cup ^\infty _{k=1} X_k $  and $Y = \cup ^\infty _{k=1} Y_k $    where
$X_1 \subset X_2 \subset X_3 \subset \ldots $ and
$Y_1 \subset Y_2 \subset Y_3 \subset \ldots $ are closed. Let $a \in X$ and $c \in Y $ be
arbitrarily  chosen.  Suppose that $X \cup Y \supset ] a,c [$. Our assumptions imply that 
$ X \cap \, ] a,c [$ and $ Y \cap \, ] a,c [$
are nonempty and that every $ x \in X $  is a
bilateral accumulation point of $X$ and every $y \in Y $ 
is a bilateral accumulation point of $Y$. We can
construct by induction an infinite sequence of nonnegative 
integers $ 0=k_0<k_1 < k_2  < \ldots $  and  the  sequences
$\{ a_r \}_{r=0}^\infty $ and $\{ c_r \}_{r=0}^{\infty} $
of real numbers such that $a_0=a, c_0=c$ and,  for  every  positive  integer  $r$, 

$$
\ \   a_{r} \in X_{{k}_{r}} \cap \, ] \, a_{r-1},c_{r-1} \, [  \ , \  c_{r} \in Y_{{k}_{r}} \cap \, ] \, 
a_{r-1},c_{r-1} \, [  \quad \text{and} \quad
  (X_{{k}_{r}} \cup Y_{{k}_{r}})   \cap   \left] \,  a_{r},c_{r} \, \right[ \, = \, \emptyset .
$$
as follows. Assume that $a_0=a, c_0=c, k_0=0$ and that $a_{r-1}, c_{r-1}, k_{r-1}$ (for a positive integer $r$) have been constructed. Choose
$\tilde a_{r} \in X \cap \, ] \, a_{r-1},c_{r-1} \, [$ and $\tilde c_{r} \in Y \cap \, ] \, a_{r-1},c_{r-1} \, [$
arbitrarily and an integer $k_r$ so large that $k_r>k_{r-1}$, $\tilde a_{r} \in X_{k_r}$ and $\tilde c_{r} \in Y_{k_r}$. As $[\tilde a_{r}, \tilde c_{r}] \cap X_{k_r}$ and $[\tilde a_{r}, \tilde c_{r}] \cap Y_{k_r}$
are two disjoint compact sets, we can choose their points $a_r$ and $c_r$, respectively, such that they realize the distance between these sets. Then we have $a_{r} \in X_{{k}_{r}} \cap \, ] \, a_{r-1},c_{r-1} \, [$, $c_{r} \in Y_{{k}_{r}} \cap \, ] \, a_{r-1},c_{r-1} \, [$ and
$(X_{{k}_{r}} \cup Y_{{k}_{r}})   \cap   \left] \,  a_{r},c_{r} \, \right[ \, = \, \emptyset$.

Now it is easy to see that for our constructed sequence of intervals $[a_r, c_r]$ we have
$$ \emptyset \not \, = \, \bigcap_{r=1}^\infty \, \left[ \, a_r,c_r \, \right] \, = \bigcap_{r=1}^{\infty} \, \left] \, a_{r},c_{r} \, \right[ \, \subset \, ] \, a,c \, [ \, \setminus %
    \, ( X \cup Y ) \ . $$
That completes the proof.

\end{proof}

{\bf 4.4 Definition.}\
\rm As  the density of the ball $ B(0,1) \subset \rmoc n  $  
is equal to  $\frac 12$ at every point of its boundary,
we can fix for any positive integer  $k$ a constant $\delta(k)$
(depending only on $k$ and on dimension $n$) such that
$ 0 < \delta(k) \leq 1$  and
$$\m n[B(e_{1},\delta(k)) \cap B(0,1)] \geq \frac{V(n)}{2} \left( 1 - \frac{1}{8k}\right)\left[\delta(k)\right]^{n}  .$$
As the function
$$ y  \longmapsto \m n[B(y,\delta (k)) \cap B(0,1)] $$
is continuous on  \rmoc n  we can fix for $k$ and  $\delta(k)$ as above
a constant $ \varepsilon(k) > 0 $  such that
$$ \m n[B(y,\delta (k)) \cap B(0,1)] \geq \frac{V(n)}{2} \left( 1 - \frac{1}{4k}\right)[\delta(k)]^{n} $$
whenever $ y \in [e_{1},(1 + \varepsilon(k)) e_{1} ]$.

According to the homogeneity and the invariance under
Euclidean isometries of \m n we see that
$$ \m n[B(y,\delta (k)) \cap B(x,r)] \geq \frac{V(n)}{2} \left( 1 - \frac{1}{4k}\right)[\delta(k)r]^{n} $$
whenever $k$ is a positive integer, $0 < r < \infty , x \in \rmoc n  $ and
$y \in B(x,(1 + \varepsilon(k))r) \setminus U(x,r) $.

\bigskip

{\bf 4.5 Theorem.}\
 Let \rpodm \Omega   be nonempty open, \rpodm A be  arbitrary
 and $\tau \in S^{n-1} $. Then
 $$ \pda = \mitau e' = \mitau pr'  . $$

\begin{proof}
 Since $\fr pr' \subset \fr e'$, it is sufficient to prove the inequalities
$$ \mitau e' \leq \pda \leq \mitau pr' . $$
\begin{itemize}
 \item[(i)]
 To prove the first inequality we may assume  $ \pda < \infty$
  and therefore we may assume  $A \cap \Omega$  is  \m n  measurable.
  According to Lemma 4.2(c) we have
  $$  \en \leq \mad \quad \text{for} \quad H_{n-1} \quad \text{a.e.} \quad z \in \rnt  . $$
  By using of Lemma  3.4  and  3.6  we see after
  the integration of the above inequality that

  $$
   \quad \, \ \    \mitau e'   =  \int\limits_{\rnt} \en \, dH_{n-1}(z)
    \leq  \int\limits_{\rnt} \mad(z) \, dH_{n-1}(z) = \pda .
  $$                                                                  
 \item[(ii)]
  To prove the inequality $$\pda \leq \mitau pr'$$ it is sufficient to
  assume that $\tau = e_n$ and
    $$\mien pr' < \infty  . $$
  We see that then  $\m n (\drr pr') = 0 $,
  according to  3.2 the set $ \Omega \cap A $ is \m n
  measurable and   $A \cap \Omega \sim \Omega \cap \ip $. Hence
      $$\pez n'(A) = \pez n'(\ip)$$
  and  according  to   3.4--3.8
  it is sufficient to prove that
 \begin{align}
  \madn (z) \leq \enpr \quad \text{for} \quad H_{n-1} \quad \text{a.e.} \quad z \in \real^{n-1}(n). 
  \end{align}
  \end{itemize}

 For any positive integer $k$ we put
 \begin{align*}
   A(k) & = \left\{ \ x \in \rmoc n   \: \m n (B(x,r) \setminus A ) \leq
           \frac{V(n)}{2}\left(1-\frac{1}{k}\right)r^{n} \quad \text{if} \quad r \in \left]0,\frac 1k \right] \ \right\},\\
   C(k) & = \left\{ \ x \in \rmoc n  \: \m n (B(x,r) \cap A ) \leq
           \frac{V(n)}{2}\left(1-\frac{1}{k}\right)r^{n} \quad
                         \text{if} \quad r \in \left]0,\frac 1k \right] \ \right\}.\\
  \end{align*}
 Obviously  $A(k)$  and $C(k)$  are closed and $A(k) \uparrow \ip,
 C(k) \uparrow int_{pr} A^c $  with  $k \uparrow + \infty $.
 For any pair of positive integers $(k,m)$ we put
 \begin{align*}
  A^+(k,m) & =  \left\{ \ x \in A(k)  \colon  Q\left(x,\frac 8m \right) \subset \Omega,
              \left] x,x + \frac 8m e_n \right[ \subset int_{pr} A^c \ \right\},\\
  A^-(k,m) & =  \left\{ \ x \in A(k) \colon Q\left(x,\frac 8m \right) \subset \Omega  ,
             \left]x,x - \frac 8m e_n \right[ \subset int_{pr}A^c \ \right\},\\
  C^+(k,m) & =  \left\{ \ x \in C(k) \colon Q\left(x,\frac{8}{m}\right) \subset \Omega  ,
             \left]x,x + \frac 8m e_n \right[ \subset \ip \ \right\},\\
  C^-(k,m) & =  \left\{ \ x \in C(k) \colon Q\left(x,\frac{8}{m}\right) \subset \Omega  ,
            \left]x,x - \frac 8m e_n \right[ \subset \ip \ \right\},\\
  B & =   \bigcup_{k=1}^\infty \bigcup_{m=1}^{\infty}
       \left[A^+(k,m) \cup  A^-(k,m) \cup  C^{+}(k,m) \cup  C^{-}(k,m) \right].
 \end{align*}
 To prove (1) it is sufficient to prove that
 \begin{align}
    \madn (z)  \leq & \,\enpr \quad
                             \text{if }  z \in \real^{n-1}(n) \setminus p_{n}(B) , \, \text{and}  \\
  H_{n-1}[\, p_{n}(B) \, ]   = & 0 \ . 
 \end{align}
 
 Firstly  we  make  the  following  observation:
 
 \noindent If $z \in \real^{n-1}(n) \setminus p_n(B) $ then the assumptions of 4.3 are
 fulfiled with  $L_n(z) , L_n(z) \cap \Omega  \cap \ip $ and
 $L_n(z) \cap \Omega  \cap int_{pr} A^c $ instead of $\real \, , \, X \text{ and }  Y , $
 respectively. Therefore for such $z$ there exists
 $ b \in \ ]a,c[ \ \cap \ \fr pr' $ whenever  $a \in L_n(z) \cap \ip $ and
 $c \in L_n(z) \cap int_{pr} A^c $ are such that  $[a,c] \subset \Omega $.

 To prove~(2) we fix  $z \in \real^{n-1}(n) \setminus p_n(B) $. We may
 assume that
 \begin{align}
  \enpr < \infty 
 \end{align}
 Then even the inclusion
 \begin{align}
  M^{\ip}_{\Omega ,n} (z) \subset \fr pr' 
 \end{align}
holds. To prove it we fix a point
$  x \in (L_n(z) \cap \Omega ) \setminus \fr pr'$.
According to~(4) we may fix $ \varepsilon > 0 $ such that
$$
  \left[ x - \varepsilon e_{n}, x + \varepsilon e_{n} \right]  \subset   \Omega  \quad \text{and} \quad
  \left[ x - \varepsilon e_{n}, x + \varepsilon e_{n} \right]  \cap  \fr pr'  =  0  .
$$
According to the observation made above, we get either
$$
  \left[ x - \varepsilon e_{n}, x + \varepsilon e_{n} \right]  \subset  \ip \quad \text{or} \quad
  \left[ x - \varepsilon e_{n}, x + \varepsilon e_{n} \right] \subset  int_{pr} A^c  \subset   (\ip)^c.
 $$
 Both cases imply that $x$ does not belong to $M^{\ip}_{\Omega ,n}(z)$.
 This completes the proof of~(5)
 and (2).
 
 To prove~(3) we fix the positive integers  $k $, $m$  and we
 shall prove that
\begin{align}
 H_{n-1}\{ \ p_{n}[A^{+}(k,m)] \ \} = 0 \ . 
\end{align}
(One can analogously prove that $p_{n}[A^{-}(k,m)], p_{n}[C^{+}(k,m)]$ and
$p_{n}[C^{-}(k,m)]$ are $H_{n-1}$ null sets.)

 To prove~(6) we put for any integer $s$
 $$ A^+(k,m,s) = \left\{ \ x \in A^+(k,m) \, \colon \, \frac{s-1}{m} < x_n \leq \frac sm  \ \right\} $$
 and assume, on the contrary, that for some fixed $s$
 we have
 \begin{align}
  H_{n-1}\{ \ p_{n}[\apl] \ \} > 0  . 
 \end{align}
 From Lebesgue outer density theorem we can  fix $z_{0} \in p_{n}[A^{+}(k,m,s)] $
 which is an outer density point (in the space $\real^{n-1}(n)$ ) of $p_{n}[A^{+}(k,m,s)] $.
 
 For every $z \in p_{n}[\apl] $  obviously there exists the
 point $ x \in \apl $  (depending on $z$)  such that
 $$
  \quad \, \ \  p_n(x) = z ,\ \, \, Q\left(x,\frac 8m \right)  \subset \Omega \quad \text{and} \quad
  \left] \zsz , \zsm ,7' \right[ \subset \left] x,\zsm ,7' \right[  \subset  int_{pr}A^c.
 $$
 We put $ x_1 = \zsm 0,1' $. According to the choice of $z_0$
 we can fix positive $r_{0}$ such that  $r_{0} \leq \frac 1m , r_{0} \leq \frac{1}{k}  $ and
 \begin{align}
  \frac{1}{V(n-1)r_{0}^{n-1}}H_{n-1} \{ \, p_{n} [U(x_{1},r_{0})] \cap p_{n}[\apl] \, \}
 \geq  1 - \frac{V(n)}{16kV(n-1)} [\delta (k) ]^{n}  ,
 \end{align}
 where $\delta(k)$ is the constant from  4.4. Putting
  $ S = p_{n} [ U(x_{1} , r_{0} )] $,
 from (8) we get
 \begin{align}
   H_{n-1} \{ \ S \cap p_{n}[\apl] \ \} \geq H_{n-1}(S) - \frac{V(n)}{16k}[\delta(k)]^{n}r^{n-1}_{0}  . 
\end{align}
 According to the choice of $x_{1}$ and \rn \ we see that
   $ \uxr 1' \cap \apl = 0 \ . $
 We can define the number $t_{0} \in \ins $ by the formula
 $$
   t_{0} = \sup \left\{ \ t \in \ins \: U(z_{0}+te_n,r_{0}) \cap \apl \not = \emptyset \ \right\}
 $$
 and we put $x_{0} = z_{0} + t_{0}e_{n} $ . The ball \uxrn \ has the
following properties:
\begin{align}
  L_n \cap \uxrn \subset int_{pr}A^c \quad \text{whenever}
  \quad     z \in p_{n}[\apl] \ , 
 \end{align}
 $\uxrn \subset \Omega  $, especially  $ A \cap \uxrn \ $ is \m n  measurable,
 $$
  \apl \cap [ \, B(x_{0},(1+\varepsilon)\rn) \setminus \uxrn \, ] \not= \emptyset
 $$
 whenever $\varepsilon > 0$.

  We fix some $ y \in \apl \cap [ \, B(x_{0},(1+\varepsilon(k))r_{0}) \setminus \uxrn \, ]$,
  where  $\varepsilon(k)$ is as in 4.4. 
  
  From  4.4  we see that
 \begin{align}
 \m n \{ \ [ \, B(y,\delta(k)\rn) \cap B(x_{0},\rn) \,] \ \} \geq \frac{V(n)}{2}\left(1-\frac{1}{4k}%
           \right)\left[\delta(k)\rn\right]^{n} .  
  \end{align}
 We define the function
 $$ g : \real^{n-1}(n) \longrightarrow [0,2\rn] $$
 by the formula
$$ g(z) = \m 1 \{ \ [ \, L_{n}(z) \cap \uxrn \, ] \setminus int_{pr}A^c \ \} \quad , z \in \real^{n-1}(n).$$
 According to 3.2  we have
 $$ \left[ \uxrn \setminus int_{pr} A^c\right] \quad \sim \quad  \left[ \uxrn \cap A \right] \ , $$
and  by  using  of  Fubini`s  theorem we get that $g$ is
$H_{n-1}$ measurable and
\begin{align}
   \m n \left[\uxrn \cap A \right] = \int\limits_{\real^{n-1}(n)}g(z) \, dH_{n-1}(z) \ . 
\end{align}
  From (10) we see that
  \begin{align*}
 g(z)& =  0 \quad \text{whenever} \quad z \in p_n[\apl] , \quad \text{and obviously} \\
 g(z)& =  0 \quad \text{whenever} \quad z \in \real^{n-1}(n) \setminus S  .
 \end{align*}
 Especially the set
 $$ \{ \ z \in \real^{n-1}(n) \: g(z) > 0 \ \} = S \setminus \{ \ z \in S \: g(z) = 0 \ \} $$
 is $H_{n-1}$ measurable and from (9) we get
  \begin{align}
& H_{n-1}\{ \ z \in \real^{n-1}(n) \, \: \, g(z) > 0 \ \} =   H_{n-1}(S) - H_{n-1} \{ z \in S \, \: \, g(z) = 0 \} \leq  \\
& \leq  H_{n-1}(S) - H_{n-1}\{S \cap p_{n}[\apl] \} \leq \frac{V(n)}{16k}[\delta(k)]^{n}\rn^{n-1}\notag
\end{align}
From (12) and (13) we see that
\begin{align}
 \m n [\uxrn \cap A ]
 \leq  2\rn H_{n-1} \{ z \in \real^{n-1}(n) \, \: \, g(z) > 0 \} \leq \frac{V(n)}{8k}[\delta(k)\rn]^{n}  .
\end{align}
  We see that
 \begin{align}
 \m n [ \, B(y,\delta(k)\rn) \setminus A \, ]
   \geq  \m n [ \, B(x_{0},\rn) \cap B(y,\delta(k)\rn) \, ] - \m n [ \, A \cap \uxrn \, ]   .
 \end{align}
 According to (11), (14) and (15) we eventually get
 \begin{align}
 \m n  [ B(y,\delta(k)\rn) \setminus A ] \geq \frac{V(n)}{2}\left(1-\frac{1}{2k}\right)\left[\delta(k)\rn\right]^{n} \ . 
 \end{align}
As $ y \in \apl \subset A(k) $ and $\delta(k)\rn \leq \rn \leq \frac{1}{k} \, , $ the
equality (16) contradicts with our definition of~$A(k)$. Hence
the assumption made in (7) leads to the contradiction and
consequently (6) and (3) hold. This completes the proof.
\end{proof}

\bigskip
{\bf 4.6 Corollary.}\
Let $\Omega \subset \rmoc n$ be nonempty open and  $A \subset \rmoc n $
be arbitrary. Then the following are equivalent :
\begin{itemize}
\item[(i)]
   $P_\Omega(A) < \infty $.
\item[(ii)]
  There  exist  linearly   independent   vectors
   $\tau_1,\tau_2,\ldots,\tau_n \in \rmoc n $ such that
     $\mu_{\tau_i} ( \Omega \cap fr_{pr} A ) < \infty \quad \text{for } i = 1,2,\ldots,n$.
\end{itemize}

\bigskip


{\bf 4.7. Lemma.}\ Let $\Omega \subseteq \mathbb{R}^n$ a  nonempty open set and $\Phi=(\Phi_1, \dots, \Phi_n)$ be an   $\mathbb{R}^n$
valued Borel measure on $\Omega$ with finite total variation. For any $\tau=(\tau_1, \tau_2, \dots, \tau_n)\in S^{n-1}$ let
$\Phi_\tau$ stand for the signed Borel measure $\sum^n_{i=1} \tau_i\Phi_i$. Then 
$$\| \Phi\|=\frac{1}{2V(n-1)}\int_{S^{n-1}}\| \Phi_\tau\| \,dH_{n-1}(\tau).$$

\begin{proof}
Let $v\colon \Omega\to \mathbb{R}^n$ be the Radon-Nikodym derivative of $\Phi$ with respect to its variation measure $|\Phi|$. Then $v$ is a $|\Phi|$ measurable $\mathbb{R}^n$ valued function and $|v(x)|=1$ for $|\Phi|$ a.e.\ $x\in\Omega$ (see \cite[2.5.12]{FH3}). 

As $\Phi_\tau(B) =\int_B \tau \circ v(x) \,d|\Phi|(x)$ for any Borel $B\subseteq \Omega$, clearly
$$\|\Phi_\tau\| =\int_{\Omega}| \tau \circ v(x) | \,d|\Phi|(x).$$
Integrating over $S^{n-1}$ and using Fubini's theorem we get
$$\int_{S^{n-1}}\| \Phi_\tau\|\,dH_{n-1}(\tau) =\int_{\Omega}\left(\int_{S^{n-1}}|\tau \circ v(x) |\, dH_{n-1}(\tau)\right)d|\Phi|(x).$$
As $S^{n-1}$ and $H_{n-1}$ are invariant under orthonormal transformations of $\mathbb{R}^n$,
$$ \int_{S^{n-1}}|\tau \circ w |\, dH_{n-1}(\tau) =|w| \int_{S^{n-1}} |\tau_1|\, dH_{n-1}(\tau)=2V(n-1)|w|\ \ \text{for any } w\in \mathbb{R}^n.$$
(See \cite[3.2.13]{FH3} for the exact values of constants $V(n-1)$ and $\int_{S^{n-1}}|\tau_1|\, dH_{n-1}(\tau)$.) Hence
$$ \int_{S^{n-1}}\| \Phi_\tau\| \, dH_{n-1}(\tau) =2V(n-1)\int_{\Omega} |v(x)|\,d|\Phi|(x)= 2V(n-1)\|\Phi\|,$$
that completes the proof.
\end{proof}

{\bf 4.8 Theorem.}\
Let $ \Omega \subset \rmoc n $ be nonempty open and $ A \subset \rmoc n$
be arbitrary. Then
$$P_{\Omega} (A)   = \frac{1}{2V(n-1)} \int_{S^{n-1}} \pda \, dH_{n-1} (\tau ).$$ 
\begin{proof}
If $P_{\Omega} (A)  =+\infty$ then clearly $\pda =+\infty$ for $H_{n-1}$ a.e. 
$\tau \in S^{n-1}$ and the equality holds. If $P_{\Omega} (A)  < +\infty$ 
then the gradient $D\chi^A$ in the sense of distributions over $\Omega$
of the characteristic function $\chi^A$ is an $\rmoc n$
valued Borel measure over $\Omega$ with finite total variation.
As $P_{\Omega} (A)  = ||D\chi^A||(\Omega)$ and $\pda = ||\tau \circ D\chi^A||(\Omega)$,
the equality holds due to the previous lemma applied to the restriction of the vector measure $D\chi^A$ to $\Omega.$
\end{proof}

{\bf 4.9 Theorem.}
Let $ \Omega \subset \rmoc n $ be nonempty open and $ A \subset \rmoc n$
be arbitrary. Then the following equalities hold:
$$P_\Omega (A)   = \je (\Omega \cap \fr e') = \je (\Omega \cap \fr pr' ).$$ 

\begin{proof}
Integrating the equalities from Theorem 4.5 over $S^{n-1}$ and using Lemma 3.4(ii) we get

$$\frac{1}{2V(n-1)} \int_{S^{n-1}} \pda \, dH_{n-1} (\tau ) = \je (\Omega \cap \fr e') = \je (\Omega \cap \fr pr' ).$$
Due to Theorem  4.8 the first term is equal to $P_\Omega(A)$. That completes the proof.
\end{proof}

{\bf 4.10 Remark.}\ The equality $P_\Omega (A)   = \je (\Omega \cap \fr e')$ for an arbitrary set
$ A \subset \rmoc n$ is known (see \cite[4.5.6]{FH3} and \cite[4.5.11]{FH3}), but our simple and 
self-contained proof does not depend on the deep results of De Giorgi, Federer and Volpert on the 
sets with finite perimeter, and it can be of independent interest. It is easy to combine our results
with other known facts and then replace integralgeometric measure $\je$ in the theorem above by the 
Hausdorff measure $H_{n-1}$ and to prove in full generality that also
$P_\Omega(A) = H_{n-1} (\Omega \cap \fr e') = H_{n-1} (\Omega \cap \fr pr') .$

\bigskip

{\bf 4.11 Some open questions.}\ We have seen that there is a variety of notions of
`measure theoretic boundary' that play an important role in integralgeometric representations
of various notions of variation of a general set $ A \subset \rmoc n$. We demonstrated 
this here using the essential boundary, and the slightly finer preponderant boundary.
While for the sets of bounded variation there is plenty of such notions of boundary that can be used,
much less is known about which notions of 'boundary' can be used for integral representations
of variations of more general sets. Even for the usual notion of the perimeter $P(A)$ of a set
$ A \subset \rmoc n$ we aim to understand for which notions of 'fine boundary', $fr_{\it fine}(A)$, we can say that $P(A)$ is equal to $(n-1)$-dimensional measure of $fr_{\it fine}(A)$ for fully general sets $ A \subset \rmoc n$. One of natural choices for such finer notions of `boundary' that need to be understand for general sets is the following `strong boundary',
$$fr_s(A) =  \{ \ x \in \rmoc n :\ \underline{d}(x, A) > 0 \quad \text{and} \quad \underline{d}(x, A^c) > 0 \ \}.$$
Or one can suggest its finer version, $fr_{s, \delta}(A)$ for $0<\delta \leq 0.5$,
$$fr_{s, \delta}(A) =  \{ \ x \in \rmoc n :\ \underline{d}(x, A) \geq \delta \quad \text{and} \quad \underline{d}(x, A^c) \geq \delta \ \}.$$
For finite variation sets these `boundaries' can be used to represent their variation, but what can be
said about them for general sets? Is their $(n-1)$-dimensional measure always equal to $P(A)$, or can it be `small' for some set $A$ of infinite variation? We know
the answers in dimension $n=1$ only, for higher dimensions these questions about 'strong boundaries' of general sets are open.

\end{document}